\documentclass[preprint,11pt,sort&compress]{elsarticle}
\usepackage{mathrsfs}
\usepackage[top=3.1cm,bottom=3.cm,left=3.cm,right=3.cm,headsep=0.5cm]{geometry}

\usepackage{circledsteps}
\usepackage{lipsum}

\usepackage{graphicx}%
\usepackage{multirow}%
\usepackage{amsmath,amssymb,amsfonts}%
\usepackage{amsthm}%
\usepackage{mathrsfs}%
\usepackage[title]{appendix}%
\usepackage{xcolor}%
\usepackage{textcomp}%
\usepackage{manyfoot}%
\usepackage{booktabs}%
\usepackage{algorithm}%
\usepackage{algorithmicx}%
\usepackage{algpseudocode}%
\usepackage{listings}%
\usepackage{setspace}
\newtheorem{theorem}{Theorem}[section]
\newtheorem{proposition}{Proposition}[section]%
\newtheorem{lemma}{Lemma}[section]
\newtheorem{corollary}{Corollary}[section]

\newtheorem{example}{Example}[section]%
\newtheorem{remark}{Remark}[section]%
\newcommand{\Real}{\mathbb{R}}

\newcommand*{\E}{\mathbb{E}}
\newcommand{\F}{\mathscr{F}}
\newcommand{\Om}{\mathrm{\Omega}}

\newcommand{\N}{\mathcal{N}}
\newcommand{\Dpdt}{ {\rm d}\mathbb{P}\otimes {\rm d}t-\mathrm{a.e.} }
\newcommand{\rd}{{ \rm d}}

\newcommand{\Hh}{\mathscr{H}}
\newcommand{\Ss}{\mathscr{S}}

\DeclareFontFamily{U}{bbold}{}%
\DeclareFontShape{U}{bbold}{m}{n}{<-5>bbold5<6>bbold6<7>bbold7
	<8>bbold8<9>bbold9<9-11>bbold10<11-14>bbold12<14->bbold12}{}%
\DeclareSymbolFont{bbold}{U}{bbold}{m}{n}
\DeclareSymbolFontAlphabet{\mathbbold}{bbold}

\begin{document}

\begin{frontmatter}

		\title{
			 Converse comparison of $g$-expectations and uniqueness of generators of BSDEs: a novel approach and some sharp results \\
 } 
 	\author[shanshi]{Xuejun  Shi}
 	\ead{shixuejun216@126.com}
 	\author[kuangda]{Shengjun Fan\corref{cor}}

 	\ead{shengjunfan@cumt.edu.cn}

 	\cortext[cor]{Corresponding author}
	
 	\author[kuangda]{Long Jiang}
 	\ead{jianglong365@cumt.edu.cn}

 	\address[shanshi]{School of Mathematics and Statistics, Shandong Normal University, Jinan, 250358, PR China}
 	\address[kuangda]{School of Mathematics, China University of Mining and		Technology, Xuzhou, 221116, PR China}

	\begin{abstract}
This paper investigates inverse problems for backward stochastic differential equations (BSDEs), with a particular focus on the converse comparison theorem and the uniqueness of generators solely through the initial values of the corresponding solutions. To this end, we develop   a forward-SDE representation approach that combines a contradiction argument with a representation of generators in terms of solutions to associated forward stochastic differential equations (SDEs). As one of the main consequences, we prove that an ordering between two $g$-expectations necessarily implies the pointwise comparability of their generators,   without requiring any additional structural assumptions beyond the standard Lipschitz condition.  This removes the continuity, $y$-independence, and deterministic assumptions required in previous works,
thereby yielding a complete solution to the converse comparison problem for $g$-expectations under the standard Lipschitz framework and, at the same time, establishing a general uniqueness theorem for  generators of BSDEs.

	\end{abstract}

	\begin{keyword}
		Backward stochastic differential equation \sep Uniqueness  \sep Converse comparison theorem \sep Inverse problem \sep $g$-expectation
	\end{keyword}
\end{frontmatter}

\section{Introduction}

The general nonlinear form of backward stochastic differential equations (BSDEs)  was first introduced by Pardoux and Peng \cite{PardouxPeng-1034}, 
\[
Y_t=\xi +\int_t^T g(s,Y_s,Z_s)\,\rd s -\int_t^T Z_s \cdot  \rd B_s, \qquad \forall t\in [0,T],
\]
whose adapted solution is typically denoted by
$
(Y_t(g,T,\xi),Z_t(g,T,\xi))_{t\in[0,T]},
$
emphasizing the dependence on the generator $g$, terminal time $T$  and terminal condition $\xi$.
The corresponding linear equation   had been investigated earlier by Bismut \cite{bismut1973conjugate} as the adjoint equation arising from the duality approach to stochastic control.
Since their introduction,   BSDEs have become a fundamental tool in modern  stochastic analysis and have found numerous applications in mathematical finance, stochastic control, and partial differential equations, see \cite{bismut1975growth,Duffie1992,peng1992nonlinear,peng1993backward,ElKarouiPeng-50,yan1997suiji_en} and the references therein. Under the standard Lipschitz condition, BSDEs are well posed, and many fundamental aspects of the theory, including comparison principles, stability properties, regularity results,  and the associated nonlinear expectations and risk measures, have been   extensively  studied over the past three decades.

One of the fundamental results in the theory of BSDEs is the comparison theorem, which states that the ordering of terminal conditions and generators is transferred to the first components of the corresponding solutions. This result plays an important role in both theoretical development  and practical applications, such as the construction of solutions for BSDEs with non-Lipschitz generators and the analysis of stochastic control problems.
Building on this foundation, Peng  \cite{Peng-564} later introduced  a special class of nonlinear expectations,   $g$-expectations,   
 under the normalization condition   $g(t,y,0)=0$,  as follows:
\[
\mathcal{E}_g[\xi]:=Y_0(g,T,\xi), \qquad \forall \xi\in L^2(\mathscr{F}_T, \mathbb{R}).
\]
The $g$-expectation preserves most desirable properties of   classical expectations  except   linearity.
Subsequently, the notions of $g$-martingale, $g$-supermartingale  and $g$-submartingale were developed based on the $g$-expectation. Using these tools, \cite{Peng-39} established a version of Doob’s decomposition for $g$-supermartingales and   \cite{ChenChen-735} investigated their relationships with other forms of nonlinear expectation, such as Choquet expectation. Furthermore, various extensions of classical results to the nonlinear framework were obtained, including the maximal (or minimal) inequality for martingales\textsuperscript{\cite{Wang2009}}, the downward-crossing inequality\textsuperscript{\cite{Chen2000}}, and the BDG inequality for $g$-martingales\textsuperscript{\cite{Faidi2024}}.
Subsequently, BSDEs have been recognized as an elegant and powerful tool in mathematical finance for the pricing and risk measurement of contingent claims; see, for instance, \cite{ElKarouiPeng-50}, \cite{RosazzaGianin2006}  and \cite{Xu2014} for details. In particular, \cite{Delbaen2009} investigated dynamical concave utilities via BSDEs, while \cite{Ji2019} and \cite{Penner2014} studied dynamical risk measures for processes within the same framework. More recently, \cite{DiNunno2024} explored fully dynamic risk measures, focusing on horizon risk, time-consistency  and their connections with BSDEs.
A natural inverse problem arising from this comparison principle is whether the converse implication holds.
  More precisely, given two generators $g_1$ and $g_2$, does the ordering of the associated $g$-expectations imply the pointwise ordering of the generators, i.e.
\[
\mathcal{E}_{g_1}[\cdot]\geq \mathcal{E}_{g_2}[\cdot]  \qquad \overset{?}{\Longrightarrow} \qquad g_1 \geq g_2.
\]
This problem, known as the converse comparison problem  for $g$-expectations, is a prototypical inverse problem for BSDEs and has attracted considerable attention in the literature.
The reader is referred to \cite{Jiang2004unique,Jiang2004converse,Jiang2005a,Briand-Coquet2000,Shi2023} for more details.
Partial results have been obtained    under some additional structural assumptions on the generators. In particular, converse comparison results have been established when the generator $g$ is independent of $y$, deterministic, or continuous in time; see, for instance, \cite{Chen1998,Briand-Coquet2000, Coquet2001, Jiang2004unique}. 
However, these assumptions are restrictive and exclude many situations of practical relevance, especially in financial applications where generators are typically random and may lack regularity in time. 
 
 To the best of our knowledge, since the pioneering works of \cite{Chen1998,Briand-Coquet2000,Coquet2001,Jiang2004converse,Jiang2004unique,Jiang2005a}, a complete converse comparison theorem for general Lipschitz generators has remained unavailable. Existing approaches rely heavily on classical BSDE representation theorems or stopping-time   techniques, which require additional assumptions such as continuity in time, independence of the $y$-variable, or deterministic structure of the generator.
  In this paper, we develop a forward-SDE  representation approach to the converse comparison problem. The proposed framework can be viewed as bridging the two principal methodologies in the existing literature, namely generator representation theorems and stopping-time arguments. The key ingredient is a new representation formula that represents the generator through the  solutions to suitably constructed forward stochastic differential equations. This viewpoint bypasses the  limitations of previous methods and enables us to establish a converse comparison theorem for $g$-expectations and a general uniqueness theorem  of generators for BSDEs under the standard Lipschitz framework.
  More precisely, we prove that for almost every  $t\in[0,T)$,
 	\[
\phi_1(t,y,z)-\phi_2(t,y,z)= \lim_{\varepsilon\rightarrow 0^+} \frac{1}{\varepsilon}   \left[\varGamma_{\phi_1}^{t,y,z}(t +\varepsilon ) -  \varGamma_{\phi_2}^{t,y,z}(t +\varepsilon )  \right], \qquad \mathbb{P}-a.s.  
\] 
 This representation enables us to reduce the comparison of generators to a comparison of forward SDE trajectories, which can then be analyzed using stopping-time arguments. As a consequence, we establish a converse comparison theorem for $g$-expectations under the standard Lipschitz framework.
 
 The main contributions of this paper are threefold:
\begin{itemize}
	\item We establish a converse comparison theorem for $g$-expectations under the standard Lipschitz condition, removing the continuity, $y$-independence, and deterministic assumptions required in previous works.
	
	\item We introduce a novel representation of the generator in terms of solutions to forward SDEs, providing a powerful tool for studying inverse problems associated with BSDEs.
	
	\item We develop a unified framework that yields converse comparison results and uniqueness of generators, and further establishes equivalent characterizations between structural properties of $g$-expectations and the corresponding generators.
 \end{itemize}
Furthermore, we demonstrate that the normalization condition is indispensable for the converse comparison theorem by constructing a counterexample in its absence. We also establish a uniqueness theorem for generators based on stopping times, which complements existing uniqueness results in the literature.

This paper is organized as follows. Section~\ref{sec:pre} introduces the notations and preliminary lemmas. Section~\ref{sec:coretech} presents the key lemmas on progressively measurable processes and forward SDEs, which serve as the foundation for the uniqueness and converse comparison results.  We also establish the main results, including the converse comparison theorem, uniqueness of generators for BSDEs, and necessary and sufficient conditions characterizing several properties of BSDE solutions.
   Section~\ref{sec:proof} is devoted to the proofs of the main results. Finally, Section~\ref{sec:conclusion} concludes the paper.

\section{Preliminaries}\label{sec:pre}
 In this  paper,   let $T>0$ be a fixed finite time horizon, and let $(B_{t})_{t\geq0}$ be a standard $d$-dimensional Brownian motion defined on a complete probability space
$(\Omega, \mathscr{F},\mathbb{P})$. Let  $(\F_t)_{0\leq t\leq T}$  be the augmented natural filtration generated by $(B_t)_{t\ge0}$, namely,
\[\F_{t}\overset{\vartriangle}{=}\sigma\{B_{s}: 0\leq s \leq t\}\vee \N,\] where $\N$ denotes the family of all $\mathbb{P}$-null sets.
Throughout the paper,   $\mathbb{R}^d$ denotes the $d$-dimensional Euclidean space equipped with the  Euclidean norm $|\cdot|$.
For $t\in[0,T]$, let $L^2(\F_t,\mathbb{R})$ denote the totality of $\F_t$-measurable, square-integrable and real-valued random variables. The following classes of adapted processes are defined:
\begin{itemize}
	\item $\Hh^2(0,T;\mathbb{R}^d)$ is the set of $\mathbb{R}^d$-valued progressively measurable processes $(Z_t)_{t\in[0,T]}$ such that
	      \[\|Z\|_{\Hh^2 } :=\left(\E\left[\int_0^T|Z_t|^2\,  \rd t\right]\right)^{1/2}<+\infty. \]
	\item $\Ss^2(0,T;\mathbb{R})$ is the set of $\mathbb{R}$-valued continuous progressively measurable processes $(Y_t)_{t\in[0,T]}$ such that
	      \[\|Y\|_{\Ss^2 }:=\left(\E\left[ \sup_{0\leq t\leq T}|Y_t|^2\right]\right)^{1/2}<+\infty.\]
\end{itemize}
We will consider the   BSDE   of the following form
\begin{equation}\label{BSDE}
	Y_t=\xi+\int_t^T g(s,Y_s,Z_s)\,\rd s -\int_t^T Z_s \cdot  \rd B_s, \quad t\in [0,T],
\end{equation}
where the  coefficients $g$,   $T$ and $\xi$ are referred to as the generator,  terminal time and terminal data, respectively.
A solution to this BSDE is a pair of
progressively measurable processes $(Y_t,Z_t)_{t\in [0,T]}$ that satisfy  the equation \eqref{BSDE}, and $(Y_t)_{t\in [0,T]}$ admits a continuous modification.  To emphasize the dependence on the coefficients $g$,  $T$ and $\xi$, we often write the solution as   $(Y_t(g,T, \xi), Z_t(g, T, \xi))_{t\in [0,T]}$.
The generator is a mapping
\[g(t,\omega, y,z):[0,  T]\times\Om\times\Real \times\Real^{  d} \mapsto \mathbb{R},\]
which is a progressively measurable stochastic process for any $(y,z)$. Throughout this paper, we assume that the generator $g$ satisfies the following conditions. For notational simplicity, the dependence on  $\omega$ is omitted.

\begin{description}
	\item[(H1)]  Square Integrability:
		\[\E\left[\left(\int_0^{T}\left|g(t, 0, 0)\right|\,  \rd t \right)^{2}\right]<+\infty;\]

	\item [(H2)] Lipschitz Condition: There exists a constant $K>0$  such that $\Dpdt$,  for each $  y,y'\in \mathbb{R}$ and $ z,z'\in \mathbb{R}^{  d}$,
	      \[
		      |g(t,y,z)-g(t,y',z')|\leq K (|y-y'| +|z-z'|)  .
	      \]
\end{description}

\begin{lemma}[Existence and Uniqueness, Theorem 3.1 in  \cite{PardouxPeng-1034},  Theorem 2.2 in \cite{yan1997suiji_en},  Theorem 4.2.1 and 4.3.1 in \cite{zhang2017backward}]\label{th:exisuniq}
	Under  assumptions (H1)--(H2), for any  $\xi \in L^2(\mathscr{F}_T,\mathbb{R})$,  BSDE \eqref{BSDE} admits a unique   solution
	$(Y_t,Z_t)_{t\in [0,T]}$ in  $\Ss^2( 0,T;\mathbb{R}) \times\Hh^2( 0,T;\mathbb{R}^d)$.
\end{lemma}

\begin{lemma}[Comparison, Theorem 3.3 in \cite{yan1997suiji_en}, Theorem 2.2 in \cite{ElKarouiPeng-50} and Lemma 1.2 in \cite{Jiang2005a}]\label{th:strictcom}
	Let $g $ and $ \bar{g}$ satisfy assumptions (H1)--(H2), and let $ \xi, \bar{\xi} \in L^2(  \mathscr{F}_T, \mathbb{R}) $. Let $ (y_t, z_t)_{t \in [0,T]} $ and $ (\bar{y}_t, \bar{z}_t)_{t \in [0,T]} $ denote the corresponding unique solutions in  $\Ss^2( 0,T;\mathbb{R}) \times\Hh^2( 0,T;\mathbb{R}^d)$ to the following BSDEs, respectively:
 \[
		y_t = \xi + \int_t^T g(s, y_s, z_s) \, \rd s - \int_t^T z_s \cdot \rd  B_s, \quad 0 \leq t \leq T,
	\]
	\[
		\bar{y}_t = \bar{\xi} + \int_t^T \bar{g}(s, \bar{y}_s, \bar{z}_s) \, \rd s - \int_t^T \bar{z}_s \cdot \rd B_s, \quad 0 \leq t \leq T.
	\]
 	Then the following assertions hold:
	\begin{enumerate}[i)]
		\item If $\xi\geq \bar{\xi}$ and $ g(t, \bar{y}_t, \bar{z}_t) \geq \bar{g}(t, \bar{y}_t, \bar{z}_t), \Dpdt$, then  we have
		      \[
			      y_t \geq \bar{y}_t \quad \text{ for all }  t\in [0,T];
		      \]
		\item If, in addition, $\mathbb{P}(\xi > \bar{\xi}) > 0 $, then
		      \[
			      \mathbb{P}(y_{t} > \bar{y}_{t}) > 0 \quad \text{for all } t \in [0,T],
		      \]
		      and, in particular, $ y_0 > \bar{y}_0 $.
	\end{enumerate}
\end{lemma}
\begin{remark}
	The conclusion in  i) is equivalent to saying that, for any stopping time $\tau \in [0,T]$, we have $y_\tau \geq \bar{y}_\tau$. Moreover, if $\xi \geq \bar{\xi}$ and $\mathbb{P}(\xi > \bar{\xi}) > 0$, then for any stopping time $\tau \in [0,T]$, it holds that $\mathbb{P}(y_\tau > \bar{y}_\tau) > 0$.
\end{remark}
In certain contexts, the following additional assumption may be imposed:
\begin{description}
	\item[(H3)]  The generator satisfies  the strong normalization condition: $\Dpdt$, for each    $ y \in \mathbb{R}$,
		\[
			g(t,y,0)=0. 
		\]
\end{description}
Furthermore, under assumptions  (H2) and (H3), the   $g$-expectation and conditional $g$-expectation are defined via the   solution  to the BSDE \eqref{BSDE}  as
\[
	\mathcal{E}_g[\xi]=Y_0(g,T,\xi),\qquad \mathcal{E}_g[\xi|\mathscr{F}_\tau]=Y_{\tau}(g,T,\xi),   \text{ with }  0\leq \tau \leq T,
\]
originally introduced by Peng \cite{Peng-564}.
These nonlinear expectation operators have the following fundamental properties.
\begin{lemma}[Lemma 36.6 in \cite{Peng-564}, Lemma 1.1 in \cite{Briand-Coquet2000}]\label{lem:progexp}
	Let $\xi_1, \xi_2 \in L^2(\mathscr{F}_T,\mathbb{R})$, and let $\tau, \sigma \in [0,T]$ be stopping times. Then,  we have:
	\begin{enumerate}[i)]
		\item Preservation of constants:  For each constant $c$, $	\mathcal{E}_g[c] = c$;
		\item   Monotonicity:  If $\xi_1 \geq  \xi_2$, then
		      $ \mathcal{E}_g[\xi_1] \geq  \mathcal{E}_g[\xi_2] 	$;
		\item  Strict monotonicity:  If $\xi_1 \geq  \xi_2$ and $P(\xi_1 > \xi_2) > 0$, then
		      $
			      \mathcal{E}_g [\xi_1]  > \mathcal{E}_g [\xi_2]
		      $;
		\item For all stopping times $\tau \geq  \sigma$,
		      $
			      \mathcal{E}_g \left[ \mathcal{E}_g [\xi_1 | \mathscr{F}_\tau]  | \mathscr{F}_\sigma \right]
			      = \mathcal{E}_g [\xi_1 |\mathscr{F}_\sigma]
		      $;
		\item For any stopping time $\tau$,
		      $
			      \mathcal{E}_g[ \xi_1 ]
			      = \mathcal{E}_g  \left[\mathcal{E}_g [\xi_1  |  \mathscr{F}_\tau] \right].
		      $
	\end{enumerate}
\end{lemma}

\section{Main results}\label{sec:coretech}
 In this section, we present the main results and outline the core ideas,  while the corresponding proofs are deferred to the next section.   
\subsection{A Forward SDE Representation of the Generator}

 We begin with establishing the following technical result.
\begin{lemma}\label{lem:convlem}
	Let $(\phi_t)_{t\in [0,T]}$ be a $\mathbb{R}$-valued progressively measurable process
	such that \[\E\left[\int_0^T|\phi_t|\,\rd t\right]<+\infty.\]
	Then, for almost every $t \in [0,T)$, we have
	\[\begin{split}
		\phi_t
		 =L^1-\lim_{\varepsilon\rightarrow 0^+} \frac{1}{\varepsilon}\int_t^{ t+\varepsilon } \phi_s\,\rd s, 
	\end{split}\]	 
	and
	\[\begin{split}	\phi_t  = \lim_{\varepsilon\rightarrow 0^+} \frac{1}{\varepsilon}\int_t^{ t+\varepsilon } \phi_s\,\rd s, \qquad  \mathbb{P}-a.s.
	\end{split}\] 
\end{lemma}
\begin{remark}
	 The above lemma extends Lemma 4 in \cite{Shi2023} and Proposition 2.2 in \cite{Jiang2008} to a more general setting.
 \end{remark}
Let
$\phi:[0,T]\times \Omega \times \mathbb{R}\times \mathbb{R}^d \to \mathbb{R}$ satisfy   assumptions (H1)--(H2). Let $(t, y, z) \in [0,T)\times \mathbb{R}\times \mathbb{R}^d$, and consider the   forward SDE:
\begin{equation}\label{eq:forwardsde}
	\varGamma(s) = y +\int_{t}^{s} \phi\left(u,\varGamma(u),z\right)\,\rd u +\int_{t}^s z  \cdot \rd B_u,\qquad s\in [t,T].
\end{equation}
By classical results, this equation admits a unique solution in $\Ss^2( t,T;\mathbb{R})$, which we denote by $ \varGamma_{\phi}^{t,y,z}(s)$ for $s\in [t,T]$ to emphasize its dependence on the parameters $(t, y, z)$ and $\phi$.
When no confusion arises, we simply write $\varGamma(s)$  for $s\in [t,T]$.

The key idea is to replace the classical BSDE representation of generators by the following forward representation, which is more compatible with trajectory comparison arguments.
\begin{proposition}[Representation of the Generator via Forward SDEs]\label{lem:forwardrep}
	Under   assumptions (H1)--(H2) on $\phi$,   for any $(y,z)\in \mathbb{R}\times \mathbb{R}^d$ 
	we have
	\[
		\phi(t,y,z)=L^1-\lim_{\varepsilon\rightarrow 0^+} \frac{1}{\varepsilon}  \E\left[\varGamma_{\phi}^{t,y,z}(t +\varepsilon ) - y  | \mathscr{F}_t\right],
	\]
	for almost every $t\in [0,T)$.
 	Furthermore, for any $y\in\mathbb{R}$ and  almost every $t\in [0,T)$, we have
	\[ 
		 \phi(t,y,0)  =L^1-\lim_{\varepsilon\rightarrow 0^+} \frac{1}{\varepsilon}    \left[\varGamma_{\phi}^{t,y,0}(t +\varepsilon ) - y   \right]   \]
		 and
	\[	  \phi(t,y,0) = \lim_{\varepsilon\to 0^+}
		 \frac{1}{\varepsilon}
		 \left[
		 \varGamma_{\phi}^{t,y,0}(t +\varepsilon )-y
		 \right],
		 \qquad \mathbb{P}-a.s.  
   \] 
\end{proposition}
\begin{remark}
	In contrast to classical representation theorems of generators for BSDEs (see Proposition~2.3 in \cite{Briand-Coquet2000} and Lemma 2.1 in \cite{Jiang2008}), this result provides an alternative representation of the generator, expressed in terms of solutions to forward SDEs rather than solutions to BSDEs, under weaker assumptions.
 \end{remark}
By an argument similar to that for Proposition \ref{lem:forwardrep}, we obtain the following result.
 \begin{proposition}[Relations Between Solutions and Generators]\label{pro:diffsolandgen}
	Under   assumptions (H1)--(H2) on $\phi_1$ and $\phi_2$,   for any $(y,z)\in \mathbb{R}\times \mathbb{R}^d$, we have
	\[
	\phi_1(t,y,z)-\phi_2(t,y,z)= \lim_{\varepsilon\rightarrow 0^+} \frac{1}{\varepsilon}   \left[\varGamma_{\phi_1}^{t,y,z}(t +\varepsilon ) -  \varGamma_{\phi_2}^{t,y,z}(t +\varepsilon )  \right], \qquad \mathbb{P}-a.s., 
	\]
	for almost every $t\in [0,T)$.
\end{proposition} 

\begin{proposition}
	\label{lem:stoptech}
	Suppose  that  assumptions (H1)--(H2) hold for $\phi_i$, $i=1,2$. 
	Let $(y,z)\in \mathbb{R}\times \mathbb{R}^d$. If  the following statement does not hold:
	\[\phi_1(t,y,z)\geq  \phi_2(t,y,z), \quad \Dpdt,\]
	 then there exists a deterministic $t \in [0,T)$,  and a stopping time  $\tau$ with $t \leq \tau \leq T$ such that
	\[
		\varGamma_{\phi_1}^{t,y,z}(\tau) \leq  \varGamma_{\phi_2}^{t,y,z}(\tau) 
	\]
	and  
	\[
		\mathbb{P}\left(\varGamma_{\phi_1}^{t,y,z}(\tau) <  \varGamma_{\phi_2}^{t,y,z}(\tau)\right)>0.
	\]
\end{proposition}

\begin{corollary}\label{cor:genequal}
	 	Suppose  that  assumptions (H1)--(H2) hold for $\phi_i$, $i=1,2$. Let $(y,z )\in \mathbb{R}\times\mathbb{R}^d$. If the following statement does not hold:
	 \[\phi_1(t,y,z)=  \phi_2(t,y,z), \quad \Dpdt,\]
	   then there exists a $t \in [0,T)$  and a stopping time  $\tau$ with $t  \leq \tau \leq T$ such that either
	 \[
	 \varGamma_{\phi_1}^{t,y,z}(\tau) \leq  \varGamma_{\phi_2}^{t,y,z}(\tau) 
	 \text{~~and~~}
	 \mathbb{P}\left(\varGamma_{\phi_1}^{t,y,z}(\tau) <  \varGamma_{\phi_2}^{t,y,z}(\tau)\right)>0,
	 \]
	 or
	 \[
	  \varGamma_{\phi_1}^{t,y,z}(\tau) \geq  \varGamma_{\phi_2}^{t,y,z}(\tau) 
	  \text{~~and~~}
	   \mathbb{P}\left(\varGamma_{\phi_1}^{t,y,z}(\tau) > \varGamma_{\phi_2}^{t,y,z}(\tau)\right)>0.
	 \] 
\end{corollary}
\subsection{The Converse Comparison Problem for $g$-Expectations}

With the above preparations, we turn to the central problems of this study,
namely the inverse problems for BSDEs. 
We first address the converse comparison problem for $g$-expectations.
Although this question has been  solved in several special settings, it has remained open
in the general framework.
The following result shows that such a conclusion holds in full generality.
It may be viewed as a fully static version of the converse comparison theorem for BSDEs.

\begin{theorem}[Converse Comparison for $g$-expectations]\label{th:convergexp}

	Assume that the generators $g_1$ and $g_2$ satisfy   assumptions  (H2)--(H3). Then the following two statements are equivalent:
	\renewcommand{\labelenumi}{\roman{enumi})}
	\begin{enumerate}
		\item  For each $\xi  \in L^2(  \mathscr{F}_T, \mathbb{R}) $,
		      $\mathcal{E}_{g_1}[\xi] \leq   \mathcal{E}_{g_2}[\xi ]$.

		\item The generators  $g_1$ and $g_2$ satisfy  that $\Dpdt$,  $  \forall (y,z)\in \mathbb{R}\times \mathbb{R}^d$,  we have
		      \[
			      g_1(t,y,z) \leq    g_2(t,y ,z ).
		      \]

	\end{enumerate}
\end{theorem}
\begin{remark} \label{re:addcons}
	 	Assume that  conditions (H2)--(H3) hold for the generators $g_1$ and $g_2$. It is well known that the above conclusion holds if  one of the following additional conditions is satisfied:
	 \renewcommand{\labelenumi}{\alph{enumi})}
	 \renewcommand{\labelenumi}{\alph{enumi})}
	 \begin{enumerate}
	 	\item \textbf{Continuity in $t$:} The generators $ g_1 $ and $ g_2 $ are almost surely continuous in $ t \in [0,T] $ for all $ (y, z) \in \mathbb{R} \times \mathbb{R}^d $.
	 	\item \textbf{Independence of $ y $:} The generators $g_1$ and $g_2$ do not depend on $ y $, that is, $\Dpdt,$  for each  $   (y, z) \in \mathbb{R} \times \mathbb{R}^d$, we have
	 	\[
	 	  g_i(t, y, z) \equiv g_i(t, z), \quad i = 1, 2.
	 	\]
	 	\item \textbf{Deterministic generators:} Both   $ g_1 $ and $ g_2 $ are deterministic functions, namely, $\rd t -  a.e.,$ for each  $    (\omega, y, z) \in \Omega\times \mathbb{R} \times \mathbb{R}^d$,  we have
	 	\[  g_i(t, \omega, y, z) \equiv g_i(t, y, z), \quad i = 1, 2.
	 	\]
	 \end{enumerate}
	  These are the special cases of converse comparison for $g$-expectations,  see  Theorem 3.2 and Theorem 4.4 in \cite{Briand-Coquet2000}, Theorem 2.1 in \cite{Coquet2001}, and Theorem 2.1 in \cite{Jiang2004unique}.	In Theorem \ref{th:convergexp}, we remove all the additional assumptions.
	 Consequently, we obtain a complete solution to the   converse comparison problem for $g$-expectations.
\end{remark}

\begin{example}
	Let the generators $g_1$ and $g_2$  be defined by
	\[\begin{aligned}
		   g_1(t,y,z) &= [f(B_t)|y|] \wedge |z|,
		   \\
		 g_2(t,y,z) &= 2\left([f(B_t)|y|] \wedge |z|\right),
	\end{aligned}
	\]
	where $f$ is a bounded, non-negative  and discontinuous function on $\mathbb{R}^d$.
	It is evident that   $g_1\leq g_2$, and both $g_1$ and $g_2$ satisfy assumptions (H2)--(H3). Consequently, by Theorem~\ref{th:convergexp}, the inequality    $\mathcal{E}_{g_1}[\xi]\leq \mathcal{E}_{g_2}[\xi]$ for any $\xi \in L^2(\mathscr{F}_T,\mathbb{R})$  provides a characterization of the ordering $g_1\leq g_2$.
	In contrast, the generators $g_1$ and $g_2$ clearly fail to satisfy the assumptions imposed in Remark \ref{re:addcons}. They are random, not continuous in $t$ and depending on the variable $y$.
\end{example}
\begin{proposition}[]\label{pro:gexpcom}
	Under   assumptions (H2)--(H3) on the generator  $g$, we have the following assertions.
	\begin{enumerate}[i)]
	\item
	Let $c\in \mathbb{R}$. For any terminal data $\xi   \in L^2( \mathscr{F}_{T}, \mathbb{R} )$,
	it holds that \[ \mathcal{E}_g[\xi +c]\leq \mathcal{E}_{g}[\xi ]+c, \]     if and only if  the generator $g$ satisfies that $\Dpdt$,  for each  $    (y,z)\in \mathbb{R}\times \mathbb{R}^d$,
	\[  g(t,y,z)\leq g(t,y-c,z). \]
	\item
	Let $\lambda \in \mathbb{R}$ and $\lambda \neq 0$. For   any terminal data
$ \xi  \in L^2( \mathscr{F}_{T}, \mathbb{R} )$, it holds that
\[\mathcal{E}_g[\lambda\xi] \leq   \lambda \mathcal{E}_g[\xi], \]   if and only if the generator $g$ satisfies  that $\Dpdt$,  for each  $    (y,z)\in \mathbb{R}\times \mathbb{R}^d$,
\[ g(t,\lambda y,\lambda z) \leq  \lambda g(t,y,z).  \]	
	\end{enumerate}
\end{proposition}

\begin{remark}
	The conclusions of Proposition \ref{pro:gexpcom} continue to hold when the inequality sign 
	``$\leq$" is replaced by ``$\geq$" or ``$=$".
\end{remark}
\begin{remark}
	
	The constant sub-additivity and star-shapedness of the $g$-expectation 
	correspond, respectively, to the cash sub-additivity and star-shapedness 
	of the induced risk measure; see \cite{castagnoli2022star,laeven2024law,laeven20231} 
	for further details.
	Part i) of Proposition \ref{pro:gexpcom} shows that constant sub-additivity 
	of the $g$-expectation is completely characterized by the monotonicity of 
	the generator with respect to $y$. 
	Part ii) provides an analogous characterization of sub-homogeneity. 
	Moreover, star-shapedness can be viewed as a special case of sub-homogeneity 
	when the constant $\lambda$ is restricted to $[0,1]$.
\end{remark}
Without assumption (H3), Theorem \ref{th:convergexp} generally fails for BSDEs. The following example demonstrates this phenomenon.
\begin{example}\label{ex:Tvarneqeq}
	Let
	\[\begin{aligned}
		 g_1(t,y,z)=g_1(t) & =-2\mathbbold{1}_{ [0,T/2)}(t)+2\mathbbold{1}_{[T/2,T]}(t),\\
		g_2(t,y,z) =g_2(t)& =2\mathbbold{1}_{ [0, T/2)}(t)- 2\mathbbold{1}_{[T/2,T]}(t).
	\end{aligned}\]
	Since \[\int_0^T g_i(s)\, \rd s =0, \qquad   i=1,2,\] 
	it is straightforward to verify that for any   $\xi\in L^2(\mathscr{F}_T;\mathbb{R})$,
	\[\begin{aligned}
		Y_0(g_1,T , \xi) & =\E[\xi] =   Y_0(g_2,T, \xi).
	\end{aligned} \]
	Nevertheless, none of the relations $g_1 = g_2$, $g_1\geq g_2$ or $g_2\geq g_1$ holds. Hence, neither the converse comparison principle nor the uniqueness of the generator can be deduced.
\end{example}

\subsection{Uniqueness of Generators for BSDEs}

    In this subsection, we investigate a particular case of the converse comparison theorem for BSDEs, namely, the uniqueness of  generators. Under assumption (H3), together with an additional continuity condition of the generator in time $t$, the uniqueness problem for   generators of BSDEs was first studied in \cite{Chen1998}, and subsequently extended to more general settings in \cite{Jiang2004unique, jiang2006limit}, see Theorem 3.3 in \cite{Jiang2004unique},  Theorem 3.4 in \cite{jiang2006limit} for details. 
 We now recall the following result about the uniqueness of generators for BSDEs.
 \begin{proposition}[Theorem 3.4 in \cite{jiang2006limit}]
 	 If for each $i=1,2$, the generator   $g_i$ satisfies   assumptions  (H1)--(H2),  and the following weak normalization condition:
 	 \[
 	 g_i(t,0,0)=0, \qquad \Dpdt,
 	 \]
 	 then the following two statements are equivalent:
 	 \begin{enumerate}[i)]
 	 	\item For any deterministic time $t \in [0,T]$,  and  any $  \xi  \in L^2(\Omega, \mathscr{F}_t )$,  we have
 	 	\[Y_0(g_1, t, \xi)=  Y_0(g_2, t, \xi). \]
 	 	\item The generators $g_1$ and $g_2$ satisfy  that $\Dpdt$,  for each   $  (y,z)\in \mathbb{R}\times\mathbb{R}^d$, we have
 	 	\[
 	 	g_1(t,y,z) = g_2(t,y,z).
 	 	\]
 	 \end{enumerate}
 \end{proposition}
 However, this result does not hold in general if the weak normalization condition
 is removed. To see this, consider the BSDE with the generator $g_1$ defined in Example \ref{ex:Tvarneqeq}, and define
 \[
 g_3(t,y,z):=g_1(t,y,z)+B_t.
 \]
 For any deterministic $t\in[0,T]$ and $\xi\in L^2(\mathscr{F}_t,\mathbb{R})$, let $Y_0$ and $Y_0'$ be the initial values of the solutions to the following BSDEs:
 \[ 
 	Y_0=\xi + \int_0^t [-2\mathbbold{1}_{ [0,T/2)}(u)+2\mathbbold{1}_{[T/2,T]}(u)] \,\rd u 
 	-\int_0^t Z_u \cdot \rd B_u, \]
  and 
 	\[Y'_0 =\xi + \int_0^t [-2\mathbbold{1}_{ [0,T/2)}(u)+2\mathbbold{1}_{[T/2,T]}(u) +B_u] \,\rd u 
 	-\int_0^t Z'_u \cdot \rd B_u.
  \]
 By Fubini's theorem, for any deterministic $t \in [0,T]$, we obtain
 \[\begin{aligned}
 	Y_0(g_1,t,\xi)=Y_0&=\E[\xi]+\E\left[\int_0^t [-2\mathbbold{1}_{ [0,T/2)}(u)+2\mathbbold{1}_{[T/2,T]}(u) ]\,\rd u\right]\\ &=Y_0'= Y_0(g_3,t,\xi).
 \end{aligned}
 \]
Nevertheless, $g_1 \neq g_3$. This example shows that the uniqueness result fails in general without the weak  normalization condition.
However, uniqueness can still be recovered if deterministic terminal times are replaced by stopping times. This leads to the following   uniqueness result under assumptions (H1)--(H2).
\begin{theorem}[Uniqueness of Generators]\label{th:mainunique}
	
	Let two generators $g_1, g_2$ satisfy assumptions  (H1)--(H2).
	Then the following two statements are equivalent:
	\begin{enumerate}[i)]
		\item For any stopping time $\tau \in [0,T]$,  and  any $  \xi  \in L^2(\Omega, \mathscr{F}_\tau)$,  we have
		\[Y_0(g_1, \tau, \xi)=  Y_0(g_2, \tau, \xi). \]
		\item The generators $g_1$ and $g_2$ satisfy that $\Dpdt$, for each  $  (y,z)\in \mathbb{R}\times\mathbb{R}^d$, we have
		\[
		g_1(t,y,z) = g_2(t,y,z).
		\]
	\end{enumerate}
	
\end{theorem}

\begin{remark}
	  For the $g_1$ and $g_2$  defined in Example \ref{ex:Tvarneqeq}, it can also be  checked that for any $\tau \in [0,T]$ and $\xi \in L^2(\mathscr{F}_\tau, \mathbb{R})$, \[ Y_0(g_1,\tau,\xi)\leq Y_0(g_2,\tau,\xi) .\]
	 But clearly, we do not have $g_1\leq g_2$.
	 This shows that, in contrast to the previous subsection, the equality ``$=$" in Theorem \ref{th:mainunique} cannot be weakened to the inequality ``$\leq$".
\end{remark}

\begin{proposition}[]\label{pro:sdeuniq}
		Let two generators $g_1$ and  $g_2$ satisfy assumptions  (H1)--(H2). Then  we have the following assertions.
	\begin{enumerate}[i)]

		\item \textbf{(Constant Preservation)} For any stopping time $\tau \in [0,T]$ and   $c\in \mathbb{R}$,    it holds that
		\[Y_0 (g  ,\tau,c)  =   c, \]   if and only if  the generator $g$ satisfies that $\Dpdt$,  for each $  y\in \mathbb{R}$,  we have \[g(t,y,0)=   0  .\]
		\item \textbf{(Zero Preservation)}  For any stopping time $\tau \in [0,T]$,    it holds that
		\[Y_0 (g  ,\tau,0)  =   0,  \]   if and only if   the generator $g$ satisfies  that  $ \Dpdt$,   $g(t,0,0)=  0.$ 
		\item \textbf{(Translation Invariance)}  For any stopping time $\tau \in [0,T]$, terminal data
		$ \xi  \in L^2( \mathscr{F}_{\tau}, \mathbb{R} )$ and any constant $c\in \mathbb{R}$, it holds that 
		\[Y_0 (g  ,\tau, \xi +c)  =   Y_0 (g  ,\tau, \xi) +c,\]  if and only if the generator $g$ is independent of $y$,  that is,  $\Dpdt$, for each  $  (y,z)\in \mathbb{R}\times\mathbb{R}^d$, we have
		\[  g(t,  y,  z) \equiv g(t,z).\]
				\item \textbf{(Homogeneity)}  Let $\lambda \in \mathbb{R}$ and $\lambda \neq 0$. For any stopping time $\tau \in [0,T]$, terminal data
				$ \xi  \in L^2( \mathscr{F}_{\tau}, \mathbb{R} )$,   it holds that
				\[Y_0(g,\tau,\lambda \xi) =    \lambda Y_0(g,\tau, \xi),   \]  if and only if the generator $g$ satisfies that   $\Dpdt$,  for each $  (y,z)\in \mathbb{R}\times\mathbb{R}^d$, we have
				\[ g(t,\lambda y,\lambda z) = \lambda g(t,y,z).\]
	\end{enumerate}

	\begin{remark}
  It follows from iv) of the Proposition \ref{pro:sdeuniq} that the homogeneity of the solution to the BSDE at time zero for all
	positive real numbers is equivalent to the generator $g$ being positively
	homogeneous in $(y,z)$.
	\end{remark}
			
\end{proposition}

\section{Proofs of the Main Results}\label{sec:proof}
\begin{proof}[Proof of Lemma \ref{lem:convlem}]
We first recall the following characterization of uniform integrability; see Theorem C.3 in \cite{oksendal2003stochastic} and Lemma 5 in \cite{Shi2023}. A family of functions $\{f_j\}_{j\in J}$ in  $L^1(\mathcal{M},\mathcal{A},\mu)$,   where $L^1(\mathcal{M},\mathcal{A},\mu)$ denotes the space of $\mathcal{A}$-measurable functions integrable with respect to the measure $\mu$,   is uniformly integrable if and only if there exists an increasing convex function $\psi:\mathbb{R}^+ \mapsto \mathbb{R}^+$ such that
\[
\lim_{x\rightarrow+\infty}\frac{\psi(x)}{x}=+\infty  
\quad\text{and}\quad
\sup_{j\in J}\int_{\mathcal{M}}\psi(|f_j|)\,\rd\mu<+\infty .
\]
	Since $\phi$ satisfies  $\E\left[\int_0^T|\phi_t|\,\rd t\right]<+\infty$, then
	there is an increasing and convex function $\varPhi$ with $\lim_{t\rightarrow +\infty}\frac{\varPhi(t)}{t}=+\infty$ such that
	$\E\left[\int_0^T\varPhi(|\phi_t|)\,\rd t\right]<+\infty$.
	Thus  $\E [\varPhi(|\phi_t|) ]<+\infty$, for a.e.  $t \in [0,T]$. By the Lebesgue differentiation theorem, we know that
	\[
	\limsup_{\varepsilon\rightarrow 0^+} \frac{1}{\varepsilon}\int_t^{t+\varepsilon} \E[\varPhi(| \phi_s|)]\,\rd s \leq  \lim_{\varepsilon\rightarrow 0^+} \frac{1}{\varepsilon}\int_t^{t+\varepsilon} \E[\varPhi(|\phi_s|)]\,\rd s = \E[\varPhi(|\phi_t|)]<+\infty
	\]
	holds for almost every $t\in [0,T)$.
	Also by the Lebesgue differentiation theorem we have, there exist $\bar{\Omega}\subset \Omega$ with $\mathbb{P}(\Omega\setminus\bar{\Omega})=0$, for any $\omega\in \bar{\Omega}$
	\[
	\lim_{\varepsilon\rightarrow 0^+} \frac{1}{\varepsilon}\int_t^{t+\varepsilon}  \phi_s\,\rd s= 	\lim_{\varepsilon\rightarrow 0^+} \frac{1}{\varepsilon}\int_t^{t+\varepsilon} \phi_s\,\rd s =\phi_t, \qquad \rd t- a.e.
	\]
	Thus, by Fubini's theorem, there exist a subset $S \subset [0,T)$ such that the Lebesgue measure  of $[0,T]\setminus S$ equals 0, and for each $t\in S$ we have
	\begin{equation}\label{eq:asconv}
		\lim_{\varepsilon\rightarrow 0^+} \frac{1}{\varepsilon}\int_t^{t+\varepsilon}  \phi_s\,\rd s=  \lim_{\varepsilon\rightarrow 0^+} \frac{1}{\varepsilon}\int_t^{ t+\varepsilon  } \phi_s\,\rd s =\phi_t, \qquad \mathbb{P}-a.s. 
	\end{equation}
	and
	\[
	\limsup_{\varepsilon\rightarrow 0^+} \frac{1}{\varepsilon}\int_t^{t+\varepsilon} \E[\varPhi(| \phi_s|)]\,\rd s \leq  \E[\varPhi(|\phi_t|)]<+\infty.
	\]
	For any $t\in S$, we know that there exist a constant $\delta_t>0$ such that
	\[
 \forall \varepsilon \in (0,\delta_t], \qquad	\frac{1}{\varepsilon}\int_t^{t+\varepsilon} \E[\varPhi(| \phi_s|)]\,\rd s \leq \E[\varPhi(|\phi_t|)]+1.
	\]
	By Jensen's inequality, we have for any $\varepsilon \in (0,\delta_t]$
	\[
	\begin{aligned}
		\E\left[ \varPhi\left(\left|\frac{1}{\varepsilon}\int_t^{ t+\varepsilon } \phi_s\,\rd s \right| \right) \right]
		\leq &
		\E\left[ \frac{1}{\varepsilon}\int_t^{t+\varepsilon} \varPhi\left(\left| \phi_s\right| \right)\,\rd s  \right]      \\
		=    & \frac{1}{\varepsilon}\int_t^{t+\varepsilon} \E\left[ \varPhi\left(\left| \phi_s\right| \right)\right]\,\rd s \\
		\leq & \E\left[ \varPhi\left(\left|\phi_t\right|] +1  \right)\right]<+\infty.
	\end{aligned}\]
	Thus,  $\left\{ \left|\frac{1}{\varepsilon}\int_t^{ t+\varepsilon } \phi_s\,\rd s \right|| \varepsilon \in (0,\delta_t]  \right\} $ are uniformly integrable.
	Then, combining this conclusion and \eqref{eq:asconv} we have
	\[\begin{split}
		\phi_t = L^1-\lim_{\varepsilon\rightarrow 0^+} \frac{1}{\varepsilon}\int_t^{ t+\varepsilon} \phi_s\,\rd s, \quad t \in S.
	\end{split}\]
	Furthermore, by Jensen's inequality, we can get
	\[\begin{split}
		\phi_t=L^1-\lim_{\varepsilon\rightarrow 0^+}\E\left[\frac{1}{\varepsilon}\int_t^{ t+\varepsilon} \phi_s\,\rd s\big|\mathscr{F}_t\right], \quad t \in S.
	\end{split}\]
	This completes the proof.
\end{proof}

\begin{proof}[Proof of Proposition \ref{lem:forwardrep}]
	Let $(y,z)\in \mathbb{R}\times \mathbb{R}^d$ and   $t \in[0,T)$,  and consider the following SDE
	\[
	\varGamma(s) = y  +\int_{t }^{s} \phi(u,\varGamma(u),z )\,\rd u +\int_{t }^s z \cdot \rd B_u, \qquad  s \in [t ,T].
	\]
	Under  assumptions (H1) and  (H2),   we have the drift term $\phi(t,y,z)$ is  Lipschitz continuous in $y$ and
	\[
	\E\left[\left(\int_{0}^T |\phi(t,0,z )|\,\rd t\right)^2\right]<+\infty, \qquad i=1,2.
	\]
	Thus the classical theory of SDEs ensures the existence and uniqueness of a strong solution, denoted by $\varGamma(t)$ for simplicity.
	
	For all $s\in[t ,T]$, we have
	\[
	\varGamma(s ) = y  +\int_{t }^{s} \phi(u,\varGamma(u),z )\,\rd u +\int_{t }^{s} z \cdot \rd B_u.
	\]
	Equivalently, it can be written in indicator form as
	\[
	\varGamma(s ) = y  +\int_{t }^{s}  \phi(u,\varGamma(u),z )\,\rd u +\int_{t }^{s }   z \cdot \rd B_u, \quad    s \in [t ,T].
	\]
	Thus we may replace $\varGamma(u)$ by $\varGamma(u)$ inside the integrand without changing the value:
	\[
	\varGamma(s ) - y =  \int_{t }^{s}  \phi(u,\varGamma(u)-y +y  ,z )\,\rd u +\int_{t }^{s}   z \cdot \rd B_u, \quad    s \in [t ,T].
	\]
	This implies that $\varGamma'(s):=\varGamma(s ) - y $ is the solution of the following SDE
	\[
	\varGamma'(s)=\int_{t }^s  \phi(u,\varGamma'(u)+y  ,z )\,\rd u +\int_{t }^{s }   z \cdot \rd B_u, \quad    s \in [t ,T].
	\]
	 By the standard $L^2$-estimate   for SDEs (see Theorem 1.25 in \cite{Yong2022} and Theorem 3.2.2 in \cite{zhang2017backward}), we know  that, for any $\varepsilon>0$, the following estimate holds on the interval $[t ,t +\varepsilon]$:
	\begin{equation}\label{eq:sdeesti}
		\mathbb{E}\left[ \max_{t  \leq u \leq t +\varepsilon} |\varGamma'(u)|^2 \right] \leq C \mathbb{E}\left[   \left( \int_{t }^{t +\varepsilon} |\phi(u,y ,z )| \rd u \right)^2 +   \int_{t }^{t +\varepsilon} |z|^2 \rd t   \right],
	\end{equation}
	where $C>0$ is a universal constant, allowed to change from line to line.
	It is evident that
	\begin{equation}\label{eq:yconver}
		\begin{aligned}
			  \frac{1}{\varepsilon}\E[\varGamma'(t +\varepsilon)|\mathscr{F}_{t }]                                                                                      
			=   \frac{1}{\varepsilon} \E\left[\int_{t }^{t +\varepsilon}   \phi(u,\varGamma'(u)+y  ,z )\,\rd u |\mathscr{F}_{t }\right].
		\end{aligned}
	\end{equation}
	Moreover, from (H2) and the estimate \eqref{eq:sdeesti}, it follows that
	\begin{equation}\label{eq:sde1st}
		\begin{aligned}
			 \lim_{\varepsilon\rightarrow 0}\E &\left[\left|\frac{1}{\varepsilon}\E\left[\int_{t }^{t +\varepsilon}   [\phi(u,\varGamma'(u)+y  ,z )-\phi(u, y  ,z )]\,\rd u \mid \mathscr{F}_{t }\right]
			\right|^2\right]                                                                                                                                                                                                                        \\
			\leq & C\lim_{\varepsilon\rightarrow 0} \E\left[\sup_{t \leq u\leq t +\varepsilon}|\varGamma'(u)|^2\right]=0.
		\end{aligned}
	\end{equation}
	By Lemma \ref{lem:convlem}, we have for almost every $t \in [0,T)$
	\begin{equation}\label{eq:sde2nd}
		L^1-  \lim_{\varepsilon\rightarrow 0}  \frac{1}{\varepsilon} \E\left[\int_{t }^{t +\varepsilon}   \phi(u, y  ,z )\,\rd u \mid \mathscr{F}_{t }\right]
		= \phi(t ,y ,z ).
	\end{equation}
	Thus, by combining \eqref{eq:sdeesti}, \eqref{eq:sde1st}, and \eqref{eq:sde2nd}, and passing to the limit in \eqref{eq:yconver}, we obtain, for almost every $t \in[0,T)$,
	\begin{equation}\label{eq:phirepre}
		\begin{aligned}
			 L^1-& \lim_{\varepsilon\rightarrow 0^+}\frac{1}{\varepsilon} \E[\varGamma((t +\varepsilon ) ) - y | \mathscr{F}_{t }] \\
			&=  L^1-\lim_{\varepsilon\rightarrow 0^+}\frac{1}{\varepsilon} \E[\varGamma'(t +\varepsilon )    | \mathscr{F}_{t }]   = \phi(t ,y ,z ).
	\end{aligned}\end{equation}

 Now, suppose that $z=0$.  For any $\varepsilon>0$ such that $t+\varepsilon\leq T$, we have  
\[
\varGamma(t+\varepsilon)-y =\int_t^{t+\varepsilon}\phi\big(u,\varGamma(u),0\big)\,\rd u.
\]
Adding and subtracting $\phi(u,y,0)$ inside the integrand gives
\begin{equation}\label{eq:gamaepsilon}
\begin{aligned}
	\frac{1}{\varepsilon}[\varGamma(t+\varepsilon)-y]
	= &\frac{1}{\varepsilon}\int_t^{t+\varepsilon} \big[\phi\big(u,\varGamma(u),0\big)-\phi(u,y,0)\big] \,\rd u
	+\frac{1}{\varepsilon}\int_t^{t+\varepsilon} \phi(u,y,0)\,\rd u.
\end{aligned}
\end{equation}
By the Lipschitz condition of $\phi$ in its $y$-variable, with Lipschitz constant $L$, we obtain
\[
\begin{aligned}
	 \left|\frac{1}{\varepsilon}\int_t^{t+\varepsilon} \big[\phi\big(u,\varGamma(u),0\big)-\phi(u,y,0)\big] \,\rd u\right|  
	 \leq L\sup_{t\leq u\leq t+\varepsilon} \big|\varGamma(u)-y\big|.
\end{aligned}
\]
Since $\varGamma$ has continuous trajectories and $\varGamma(t)=y$, it follows that
\[
\sup_{t\leq u\leq t+\varepsilon} \big|\varGamma(u)-y\big|\longrightarrow 0, \qquad \mathbb{P}-a.s.,
\]
as $\varepsilon\to0^+$. Hence,
\begin{equation}\label{eq:firstcon}
\lim_{\varepsilon\to0^+} \frac{1}{\varepsilon}\int_t^{t+\varepsilon} \big[\phi\big(u,\varGamma(u),0\big)-\phi(u,y,0)\big] \,\rd u=0, \qquad \mathbb{P}-a.s.
\end{equation}
Moreover, estimate \eqref{eq:sdeesti} yields
\begin{equation}\label{eq:firstconl1}
L^1- \lim_{\varepsilon\to0^+} \frac{1}{\varepsilon}\int_t^{t+\varepsilon} \big[\phi\big(u,\varGamma(u),0\big)-\phi(u,y,0)\big] \,\rd u=0.
\end{equation}	
On the other hand, the Lebesgue differentiation theorem implies that
\begin{equation}\label{eq:secondcon}
\lim_{\varepsilon\to0^+} \frac{1}{\varepsilon}\int_t^{t+\varepsilon} \phi(u,y,0)\,\rd u =\phi(t,y,0), \quad \mathbb{P}-a.s.,
\end{equation}
for almost every $t\in[0,T)$. 

Now, combining \eqref{eq:gamaepsilon},   \eqref{eq:firstcon} and \eqref{eq:secondcon} , we conclude that
\[
\lim_{\varepsilon\to0^+} \frac{1}{\varepsilon}[\varGamma(t+\varepsilon)-y]=\phi(t,y,0), \qquad \mathbb{P}-a.s.,
\]
for almost  every $t\in[0,T)$.  
Furthermore, combining \eqref{eq:gamaepsilon} and \eqref{eq:firstconl1}  with Lemma \ref{lem:convlem}, we obtain  
\[
 L^1-\lim_{\varepsilon\to0^+} \frac{1}{\varepsilon}[\varGamma(t+\varepsilon)-y] =\phi(t,y,0).
\]
This completes the proof.
\end{proof}

\begin{proof}[Proof of Proposition \ref{pro:diffsolandgen}] 
	For each fixed $(t,y,z)\in [0,T)\times\mathbb{R}\times\mathbb{R}^d$, consider the following SDEs:
	\[
	\varGamma_i(s)
	=y+\int_t^s \phi_i(u,\varGamma_i(u),z)\,\rd u
	+\int_t^s z \cdot \rd B_u,
	\quad s\in[t,T],\quad i=1,2.
	\]
	Subtracting the equation for $\varGamma_2$ from that for $\varGamma_1$, we obtain
	\[
	\varGamma_1(s)-\varGamma_2(s)
	=\int_t^s
	\big[\phi_1(u,\varGamma_1(u),z)
	-\phi_2(u,\varGamma_2(u),z)\big]\,\rd u,
	\quad s\in[t,T].
	\]
	Consequently, for every $\varepsilon\in (0,T-t]$,
	\begin{equation}\label{eq:diff3}
	\begin{aligned}
		   \varGamma_1(t+\varepsilon)-\varGamma_2(t+\varepsilon) 
		 = \int_t^{t+\varepsilon}&
		 \big[\phi_1(u,\varGamma_1(u),z)
		 -\phi_1(u,y,z)\big]\,\rd u\\&
		  + \int_t^{t+\varepsilon}
		 \big[\phi_2(u,y,z)-\phi_2(u,\varGamma_2(u),z)\big]\,\rd u \\
		 &+\int_t^{t+\varepsilon}
		 \big[\phi_1(u,y,z)
		 -\phi_2(u,y,z)\big]\,\rd u.
	\end{aligned}	 
	\end{equation}
	By the Lipschitz continuity of $\phi_1$ and $\phi_2$,
	we have
	\begin{equation}\label{eq:difffirst2}
		\begin{aligned}
			& \left|\int_t^{t+\varepsilon} 
			 \big[\phi_1(u,\varGamma_1(u),z)
			 -\phi_1(u,y,z)\big]\,\rd u 
			 + \int_t^{t+\varepsilon}
			 \big[\phi_2(u,y,z)-\phi_2(u,\varGamma_2(u),z)\big]\,\rd u  \right|\\
		&	\quad \leq \varepsilon\cdot L \left[\sup_{t\leq u\leq t+\varepsilon} \left(|\varGamma_1(u)-y| +|\varGamma_2(u)-y|\right)\right], \qquad \mathbb{P}-a.s.
		\end{aligned} 
	\end{equation}
	Since $\varGamma_1(t)=\varGamma_2(t)=y$ and
	both processes have continuous sample paths,
	the  estimate \eqref{eq:difffirst2}
	implies that its left-hand side, divided by
	$\varepsilon$, converges to zero almost surely
	as $\varepsilon\to0^+$.
	Thus, by \eqref{eq:diff3} and the Lebesgue
	differentiation theorem, we obtain
	
	\[
	\phi_1(t,y,z)-\phi_2(t,y,z)
	=\lim_{\varepsilon\to 0^+}\frac1\varepsilon\left[\varGamma_{1}(t+\varepsilon)-\varGamma_{2} (t+\varepsilon)\right],
	\qquad\mathbb{P}-a.s.
	\]
	for almost every $t\in[0,T)$. This completes the proof. 
\end{proof}

\begin{proof}[Proof of Proposition \ref{lem:stoptech}]

	According to Proposition \ref{pro:diffsolandgen}, if the inequality
	\[
	\phi_1(t,y,z)\geq   \phi_2(t,y,z),\qquad   \Dpdt
	\]
	fails to hold,  then  there  exist $t\in[0,T)$, $\delta>0$, and
	an event $\Omega_t^\delta\in\mathscr{F}_T$, with
	$\mathbb{P}(\Omega_t^\delta)>0$, such that 
	\[
	   \phi_2(t,y,z)-\phi_1(t,y,z)=\lim_{\varepsilon\to0^+}\frac{1}{\varepsilon}
	\left[
	\varGamma_{\phi_2}^{t,y,z}(t+\varepsilon)
	-
	\varGamma_{\phi_1}^{t,y,z}(t+\varepsilon)
	\right]>\delta, \quad  \text{on }\Omega_t^\delta.
 	\]
	
	For each $\omega\in\Omega_t^\delta$, the above limit implies that there exists a
	number $\eta_\omega \in (0,T-t)$, depending on $\omega$, such that
	\begin{equation}\label{eq:localstrict}
	\varGamma_{\phi_2}^{t,y,z}(s,\omega)
	-\varGamma_{\phi_1}^{t,y,z}(s,\omega)
	>\delta(s-t), \quad \text{for all } s\in (t,t+\eta_\omega]. 
	\end{equation}
	At $s=t$, however, $\varGamma_{\phi_i}^{t,y,z}(t)=y$ for $i=1,2$, and hence
	\[
	\varGamma_{\phi_2}^{t,y,z}(t)
	-\varGamma_{\phi_1}^{t,y,z}(t)=y-y=0.
	\]
	Consequently, for each $\omega\in\Omega_t^\delta$ and 
	 $s\in[t, t+\eta_\omega]$, 
	\[
	\varGamma_{\phi_2}^{t,y,z}(s,\omega)
	\geq
	\varGamma_{\phi_1}^{t,y,z}(s,\omega)+\delta(s-t).
	\]
	
	For the fixed triple $(t,y,z)\in[0,T)\times\mathbb{R}\times
	\mathbb{R}^d$, define
	\[
	\tau_\delta
	:=\inf\left\{
	s\geq t:
	 \varGamma_{\phi_2}^{t,y,z}(s)<\varGamma_{\phi_1}^{t,y,z}(s)+\delta(s-t)
	\right\}\wedge T.
	\]
   By \eqref{eq:localstrict} and the fact $\mathbb{P}(\Omega_t^\delta)>0$, we have
	\[
	\tau_\delta\geq t 
	\qquad \text{and}\qquad
	\mathbb{P}(\tau_\delta>t)\geq \mathbb{P}(\Omega_t^\delta)>0.
	\]
   Furthermore, by the definition of $\tau_\delta$ and the continuity of the
	sample paths of $\varGamma_{\phi_i}^{t,y,z}$, $i=1,2$,  it follows that
	\[
	\varGamma_{\phi_1}^{t,y,z}(s)+\delta(s-t)\leq\varGamma_{\phi_2}^{t,y,z}(s),
	\qquad s\in[t,\tau_\delta].
	\]
	In particular,
	\[
	\varGamma_{\phi_1}^{t,y,z}(\tau_\delta)\leq\varGamma_{\phi_2}^{t,y,z}(\tau_\delta).
	\]
	Moreover, on the event $\{\tau_\delta>t\}$,
	\[
	\varGamma_{\phi_1}^{t,y,z}(\tau_\delta)
	<
	\varGamma_{\phi_1}^{t,y,z}(\tau_\delta)+\delta(\tau_\delta-t)
	\leq\varGamma_{\phi_2}^{t,y,z}(\tau_\delta).
	\]
	Thus,  $ \tau_\delta$ is  the desired stopping time, and the proof
	is complete.
\end{proof}

\begin{proof}[Proof of Theorem \ref{th:convergexp}]
	The implication  ii) $\Longrightarrow$  i) follows directly from the classical comparison theorem.
	Therefore, it remains to prove i) $\Longrightarrow$ ii).
	Let $(t ,y ,z )\in [0,T)\times \mathbb{R} \times \mathbb{R}^d $, and consider the processes
	$(Y _i(s))_{s \in [t ,T]}$ for $i=1,2$ defined as the unique solution to   the following   forward SDEs:
	\begin{equation}\label{eq:eq:uniquesde1}
		Y_i(s) = y  -\int_{t }^{s} g_i(u,Y_i(u),z )\,\rd u +\int_{t }^s z \cdot \rd B_u, \qquad  s\in[t,T].
	\end{equation}
	We establish the implication by contradiction. Suppose that the condition  ii) does not hold. Then, by applying Proposition \ref{lem:stoptech} with $\phi_i = -g_i$ for $i = 1, 2$, it follows that there exist  $(t, y, z) \in [0,T)\times \mathbb{R} \times \mathbb{R}^d  $  and a stopping time $\tau$
	with $\tau \in [t,T]$,   such that
	\[Y_1(\tau) \leq Y_2(\tau),\]
	and moreover, on a subset $\bar{\Omega} \subset \Omega$ with $\mathbb{P}(\bar{\Omega} )>0$,   the strict inequality holds:
	\[ Y_1(\tau) <  Y_2(\tau) .\]
	
	We rewrite the  forward SDE \eqref{eq:uniquesde} in the backward form:
	\[Y_i(s)   = Y_i(T) + \int_s^T g_i\left(u, Y_i(u), z\right)\,\rd u - \int_s^T z  \cdot \rd B_u, \quad s \in [t ,T].\]
	Then the pairs of processes $(Y_i(s), z )_{s\in [t ,T]}$, $i = 1, 2$, are the unique solutions of BSDEs with parameters $(g_i, T, Y_i (T))$:
	\[
	\bar{Y}^i_s   = Y_i(T)+ \int_s^T g_i\left(u, \bar{Y}^i_u, \bar{Z}^i_u\right)\,\rd u - \int_s^T \bar{Z}^i_u \cdot \rd B_u, \qquad s\in[t,T].
	\]
	Hence, by the fundamental properties of $g$-expectation, we have
	\[
	\mathcal{E}_{g_1}\left[ Y_1 ( \tau)  \big\vert \mathscr{F}_{t } \right] = \mathcal{E}_{g_1}\left[ Y_1 (T) \big\vert \mathscr{F}_{t } \right] = y,
	\]
	and
	\[
	\mathcal{E}_{g_2}\left[ Y_2(\tau) \big\vert \mathscr{F}_{t } \right] = \mathcal{E}_{g_2}\left[ Y_2( T )\big\vert \mathscr{F}_{t } \right] = y.
	\]
	Thus, we   obtain
	\begin{equation}\label{eq:conkey2}
		\mathcal{E}_{g_1}\left[ Y_1(\tau) \right] = \mathcal{E}_{g_2}\left[ Y_2(\tau) \right] = y .
	\end{equation}
	But in light of the  fact that   $ Y_1(\tau) \leq  Y_2(\tau)$ with $\mathbb{P}\left(Y_1(\tau) < Y_2(\tau)\right) > 0$, it follows from the strict monotonicity as stated in iii) of Lemma \ref{lem:progexp} along with the condition i) that
	\[
	y  = \mathcal{E}_{g_1}\left[ Y_1(\tau) \right] < \mathcal{E}_{g_1}\left[ Y_2(\tau) \right] \leq  \mathcal{E}_{g_2}\left[ Y_2(\tau) \right] = y .
	\]
	This induces a contradiction. The proof is complete.
\end{proof}

\begin{proof}[Proof of Proposition \ref{pro:gexpcom}]
	We first prove (i).
	Let $c \in \mathbb{R} $ and $\xi \in L^2(\mathscr{F}_T,\mathbb{R})$, and consider the following two BSDEs:
	\[
	y_t =   \xi +c+ \int_t^T g(s, y_s, z_s) \, \rd s - \int_t^T z_s \cdot \rd  B_s, \quad 0 \leq t \leq T,
	\]
	and
	\[
	y'_t =   \xi + \int_t^T g(s, y'_s, z'_s) \, \rd s - \int_t^T z'_s \cdot \rd  B_s, \quad 0 \leq t \leq T.
	\]
	Thus we have
	\[
	y'_t +c=     \xi+c + \int_t^T   g(s, y'_s, z'_s) \, \rd s - \int_t^T   z'_s \cdot \rd  B_s, \quad 0 \leq t \leq T.
	\]
	Setting $g^c(t,y,z):=g(t,y-c,z)$,
	we know that $(y_t'+c)_{t\in[0,T]}$ is the solution of the BSDE with terminal data $\xi+c$ and generator $g^c$.
	Thus, we have
	\[
		\forall \xi \in L^2(\mathscr{F}_T, \mathbb{R}), \qquad   \mathcal{E}_{g^c}[\xi+c]=\mathcal{E}_g[\xi]+c.
\]
	Consequently, the statement that for each $\xi\in  L^2(\mathscr{F}_T,\mathbb{R})$,  $\mathcal{E}_g[\xi+c]\leq \mathcal{E}_g[\xi]+c$  is equivalent to the statement that  
	\begin{equation}\label{eq:conssubadd}
	\forall \xi \in L^2(\mathscr{F}_T, \mathbb{R}), \qquad  	\mathcal{E}_g[\xi+c]\leq \mathcal{E}_{g^c}[\xi+c].
	\end{equation} 
	 	Invoking Theorem \ref{th:convergexp},  the last statement \eqref{eq:conssubadd} is equivalent to    $\Dpdt$,
	\[
  \forall (y,z) \in \mathbb{R}\times\mathbb{R}^d, \qquad g(t,y,z)\leq g^c(t,y,z)=g(t,y-c,z).
	\]
	This is the desired conclusion.
	
	Then we prove ii).
	Let $\lambda \neq 0$ and $\xi \in L^2(\mathscr{F}_T,\mathbb{R})$, and consider the following two BSDEs:
	\[
	y_t = \lambda \xi + \int_t^T g(s, y_s, z_s) \, \rd s - \int_t^T z_s \cdot \rd  B_s, \quad 0 \leq t \leq T,	\]
	and
	\[
	y'_t =   \xi + \int_t^T g(s, y'_s, z'_s) \, \rd s- \int_t^T z'_s \cdot \rd  B_s, \quad 0 \leq t \leq T.
	\]
	Multiplying the second equation by $\lambda$, we obtain
	\[
	\lambda y'_t =   \lambda \xi + \int_t^T \lambda g(s, y'_s, z'_s) \, \rd s - \int_t^T \lambda z'_s \cdot \rd  B_s, \quad 0 \leq t \leq T.
	\]
	Define the rescaled generator $g^{\lambda}(t,y,z):=\lambda g(t,\frac{1}{\lambda}y,\frac{1}{\lambda}z)$.
	Then the above equation can be rewritten as
	\[
	\lambda y'_t =   \lambda \xi + \int_t^T   g^{\lambda}(s, \lambda y'_s, \lambda z'_s) \, \rd s - \int_t^T \lambda z'_s \cdot \rd  B_s, \quad 0 \leq t \leq T,
	\]
	According to the uniqueness of solutions to BSDEs, we know that 
	\[
	  \forall \xi \in L^2(\mathscr{F}_T,\mathbb{R}), \qquad \lambda\mathcal{E}_g[ \xi] =\mathcal{E}_{g^{\lambda}}[\lambda \xi].
	\]
		Consequently, the statement that for each $\xi\in  L^2(\mathscr{F}_T,\mathbb{R})$,  $\mathcal{E}_g[\lambda\xi  ]\leq \lambda \mathcal{E}_g[\xi] $  is equivalent to the statement that   
	\begin{equation}\label{eq:subhom} 
		\forall \xi \in L^2(\mathscr{F}_T, \mathbb{R}), \qquad  \mathcal{E}_g[\lambda \xi ]\leq   \mathcal{E}_{g^\lambda }[\lambda\xi ].\end{equation}
	Invoking Theorem \ref{th:convergexp},  the last statement \eqref{eq:subhom} is equivalent to    $\Dpdt$,
	\[
	 \forall (y,z) \in \mathbb{R}\times\mathbb{R}^d, \qquad g(t,y,z)\leq g^{\lambda}(t,y,z)=\lambda g\left(t,\frac{1}{\lambda} y, \frac{1}{\lambda}z\right),
	\]
	which implies that $\Dpdt$
	\[
	\forall (y,z) \in \mathbb{R}\times\mathbb{R}^d, \qquad g(t,\lambda y,\lambda z)\leq g^{\lambda}(t,\lambda y,\lambda z)=\lambda g(t,  y,  z).
	\]
 	The proof is thus complete.
\end{proof}

\begin{proof}[Proof of Theorem \ref{th:mainunique}]
	The implication ii) $\Longrightarrow$ i) is an immediate consequence of the comparison theorem or existence and uniqueness of solutions to BSDEs.
	Now, we prove i) $\Longrightarrow$ ii).
	Let $(t,y,z)\in [0,T)\times \mathbb{R}\times\mathbb{R}^d$, and consider the processes
	$(Y _i(s))_{s \in [t ,T]}$ for $i=1,2$ defined as the unique solution to   the following  forward SDEs:
	\begin{equation}\label{eq:uniquesde}
		Y_i(s) = y  -\int_{t }^{s} g_i(u,Y_i(u),z )\,\rd u +\int_{t }^s z \cdot \rd B_u, \qquad  
		s\in [t,T].
	\end{equation}
	Suppose ii) fails. Then by Corollary \ref{cor:genequal}, without loss of generality, we may assume that  there exist $(t, y, z)\in [0,T)\times \mathbb{R}\times\mathbb{R}^d$,
	and a stopping time $\tau \in [t ,T]$ such that
	\[Y_1(\tau) \leq Y_2(\tau),\]
	and moreover, on a subset $\bar{\Omega}\subset \Omega $ with $\mathbb{P}(\bar{\Omega} )>0$,   the strict inequality holds:
	\[	Y_1(\tau)  <  Y_2(\tau)  .\]

	Now  we study the solutions of the following two BSDEs, which are corresponding to the  forward SDEs \eqref{eq:uniquesde}:
	\[
	\begin{aligned}
		\bar{Y}^1_s & = Y_1(T) + \int_s^T g_1\left(u, \bar{Y}^1_u, \bar{Z}^1_u\right)\,\rd u - \int_s^T \bar{Z}^1_u \cdot \rd B_u, \quad s \in [0,T], \\
		\bar{Y}^2_s & = Y_2(T)+ \int_s^T g_2\left(u, \bar{Y}^2_u, \bar{Z}^2_u\right)\,\rd u - \int_s^T \bar{Z}^2_u \cdot \rd B_u, \quad s \in [0,T].
	\end{aligned}
	\]
	It follows from \eqref{eq:uniquesde} that $\bar{ Y}^i_s=Y_i(s)$ and $\bar{Z}_s=z$ for $s\in [t,T]$ and $i=1,2$.	Obviously we have
	\[
	Y_{t}\left(g_1, T, Y_1(T)\right) = \bar{Y}^1_{t} =Y_1(t)= y, \quad Y_{t}\left(g_2, T, Y_2(T)\right) = \bar{Y}^2_{t} =Y_2(t)= y .
	\]
	Thus, by Proposition 2.5 of  \cite{ElKarouiPeng-50},  we can get that
	\[
	\begin{aligned}
		Y_0\left(g_1, t, y\right) &
		= Y_0\left(g_1, T, Y_1(T)\right)  = Y_0\left(g_1, \tau, \bar{Y}^1_{\tau}\right);
	\end{aligned}
	\]                          	\[
	\begin{aligned}
		Y_0\left(g_2, t, y \right) &
		 = Y_0\left(g_2, T, Y_2(T)\right)  = Y_0\left(g_2, \tau, \bar{Y}^2_{\tau}\right).
	\end{aligned}
	\]
	Then, combining with the statement i) of  this theorem,    we have
	\begin{equation}\label{eq:solBSDE1}
		\begin{aligned}
			Y_0\left(g_2, \tau, \bar{Y}^1_{\tau}\right)  = Y_0\left(g_1, \tau, \bar{Y}^1_{\tau}\right)    = Y_0\left(g_1, t, y\right)
			= Y_0\left(g_2, t, y \right)                                                                                           = Y_0\left(g_2, \tau, \bar{Y}^2_{\tau}\right).
		\end{aligned}
	\end{equation}
	
	On the other hand,
	\[
	\bar{g}_2(t, y, z) = g_2(t, y, z) \mathbbold{1}_{[0, \tau]}(t), \quad \forall   (t,y, z) \in [0,T]\times  \mathbb{R} \times \mathbb{R}^d,
	\]
	where $\mathbbold{1}_{[0, \tau]}$ denotes the indicator function of $[0, \tau]$. It follows immediately that  $Y_0\left(g_2, \tau,\bar{ Y}^i_{\tau}\right) = Y_0\left(\bar{g}_2, T, \bar{Y}^i_{\tau}\right), i = 1, 2$,  by Proposition~2.5 in \cite{ElKarouiPeng-50}; see also the proof of Theorem~3.1 in \cite{Jiang2004unique}.
	 From the discussion above, we have
	\[
	\bar{Y}^1_{\tau}=Y_1(\tau)\leq Y_2(\tau) = \bar{Y}^2_{\tau} \text{ and } \mathbb{P}\left(\{\bar{Y}^1_{\tau} < \bar{Y}^2_{\tau}\}\right)   > 0.
	\]
	By applying Lemma \ref{th:strictcom} to the two BSDEs with common generator $\bar{g}_2$,  it follows that
	\begin{equation}\label{eq:solBSDE2}
		\begin{aligned}
			Y_0\left(g_2, \tau, \bar{Y}^1_{\tau}\right) & = Y_0\left(\bar{g}_2, T, \bar{Y}^1_{\tau}\right)       
			  < Y_0\left(\bar{g}_2, T, \bar{Y}^2_{\tau}\right)  =     Y_0\left(g_2, \tau, \bar{Y}^2_{\tau}\right),
		\end{aligned}
	\end{equation}
	which is a contradiction to \eqref{eq:solBSDE1}.
	The proof  is complete.
\end{proof}

\begin{proof}[Proof of Proposition \ref{pro:sdeuniq}]
	We first prove  i). The solutions of the BSDE are constant preserving if and only if $g(t,y,0)=0$ for all $y\in\mathbb{R}$. The sufficiency follows immediately from the uniqueness of solutions to the BSDE, and we now turn to the proof of the necessity.
	Let $(t,y)\in [0,T)\times \mathbb{R}$, and consider the processes
	$(Y (s))_{s \in [t ,T]}$ and $(Y' (s))_{s \in [t ,T]}$ defined as the unique solutions to   the following system of the forward SDEs:
	\begin{equation}\label{eq:constsde}
		\begin{aligned}
			Y(s) &= y  -\int_{t }^{s} g(u,Y(u),0 )\,\rd u +\int_{t }^s 0 \cdot \rd B_u, \qquad   s\in [t,T],\\
			Y'(s) &=y - \int_{t }^{s} 0\,\rd u +\int_{t }^s 0 \cdot \rd B_u, \qquad   s\in [t,T].	 
		\end{aligned}
	\end{equation}
	It is evident that $Y'(s)=y$ on the interval $[t,T]$.
	Now we prove the necessary part  by contradiction.  
	If  the statement that
	\[
	  \Dpdt, \qquad\forall y\in \mathbb{R}, \quad 
	 g(t,y,0)=   0.\]
	fails to hold, then by Corollary \ref{cor:genequal},  without   loss of generality,  we may assume  that there exists  $t\in [0,T)$, $y\in \mathbb{R}$
	and a stopping time $\tau$  on the interval $[t ,T]$ such that
	\[
	Y(\tau) > Y'(\tau)=y   \qquad \text{on a set }  \bar{\Omega}\subset \Omega \text{ with } \mathbb{P}(\bar{\Omega} )>0,
	\]
	and
	\[
	Y(\tau) \geq  Y'(\tau)=y .
	\]
	Note that the corresponding SDE with drift term $-g$ and parameters $(t ,y ,0)\in [0,T)\times \mathbb{R}\times\mathbb{R}^d$ can be rewritten in backward form as
	\[
	Y(s)=Y(\tau)+\int_{s}^{\tau}g(u,Y(u),0) \rd u -\int_{s}^{\tau} 0 \cdot \rd B_u, \qquad s\in [t ,\tau].
	\]
	Consequently,  combining the constant preserving property and the comparison theorem, we obtain
	\[
	\begin{aligned}
		y =Y_0(g,t,y)=Y_0(g,\tau, Y(\tau)) =Y_0(\bar{g} ,T, Y(\tau)) >Y_0(\bar{g},T,y) =Y_0(g,\tau,y )=y ,
	\end{aligned}
	\]
	where
	\[
	\bar g(s,y,z):=g(s,y,z)\mathbbold{1}_{[0,\tau]}(s),
	\qquad
	(s,y,z)\in[0,T]\times\mathbb R\times\mathbb R^d.
	\]
	This contradiction completes the proof of i).
	
	The proof of assertion ii) is analogous to that of assertion i).
	For assertion iii), by  replacing Theorem~\ref{th:convergexp} with Theorem \ref{th:mainunique}, replacing the inequality  sign ``$\leq$" with the equality sign ``$=$",   and  following  a similar argument as that in the proof of part i) of Proposition \ref{pro:gexpcom}, we deduce that
	\[
	\Dpdt, \qquad   \forall\, y,c \in \mathbb{R}, \forall z\in \mathbb{R}^d, \quad g(t,y,z)=g(t,y-c,z)=g(t,0,z).
	\]
	 Therefore, the desired conclusion that $g$ is independent of $y$ follows immediately.  
	 Similarly, the proof of assertion iv) follows the same idea as that used in the proof of part ii) of Proposition~\ref{pro:gexpcom}. For any fixed $\lambda\neq 0$  and any $\xi\in L^2(\mathscr{F}_\tau,\mathbb{R})$, the   homogeneity property
	 \[
	 Y_0(g,\tau,\lambda \xi)=\lambda Y_0(g,\tau,\xi)
	 \]
	 can be rewritten as
	 \begin{equation}\label{eq:unihom}
	 \forall\,\xi\in L^2(\mathscr{F}_\tau,\mathbb{R}),
	 \qquad
	 Y_0(g,\tau,\lambda \xi)
	 =
	 Y_0(g^\lambda,\tau,\lambda \xi),
	 \end{equation}
	 where
	 \[
	 g^\lambda(t,y,z):
	 =
	 \lambda g\left(t,\frac{1}{\lambda}y,\frac{1}{\lambda}z\right), \qquad   (t, y,z)\in [0,T]\times  \mathbb{R}\times \mathbb{R}^n.
	 \]
	 Therefore, by the uniqueness result in Theorem~\ref{th:mainunique} and \eqref{eq:unihom}, we obtain
	 \[
	 \Dpdt,
	 \qquad
	 g(t,\lambda y,\lambda z)=\lambda g(t,y,z),
	 \]
	 for  $(y,z)\in\mathbb{R}\times\mathbb{R}^d$.
	 This completes the proof.
\end{proof}
\section{Conclusions}\label{sec:conclusion}
In this paper, we developed a novel forward-SDE representation approach to inverse problems for BSDEs and applied it to the study of converse comparison for $g$-expectations and uniqueness of generators. This approach yields a fully general converse comparison theorem, a new uniqueness theorem for generators of BSDEs, and several sharp characterization results under the standard Lipschitz framework. We expect that the proposed methodology will also be useful for other inverse problems involving BSDEs and nonlinear expectations, including multidimensional, mean-field, and quadratic settings.

\section*{Acknowledgements}
	The first author is supported by the Natural Science Foundation of Shandong Province (Grant No. ZR2022MA079),   and the second author is supported by   National Natural Science Foundation of China (Grant No. 12171471),  Key Program of the
	National Natural Science Foundation of China(Grant No. 12631016).
 \small
\bibliographystyle{elsarticle-harv}
 \addtolength{\bibsep}{-1.2ex}
\bibliography{BSDE} 

\end{document}